\documentclass[11pt,leqno]{article}

\usepackage{amssymb,amsmath,amsthm}

\newtheorem{lem}{Lemma}
\newtheorem*{thm}{Theorem}

\def\n{\noindent}
 
\def\vp{\varepsilon}

\newcommand{\bb}[1]{\mathbb{#1}}

\def\NN{\bb N}

\def\H{{\cal H}}

\begin{document}

\title{The Operator Hilbert
 Space $OH$ and TYPE III von Neumann
Algebras}

\author{Gilles Pisier\thanks{Partially supported by   NSF      and
   Texas Advanced Research
 Program 010366-163}\\
Texas A\&M University\\
College Station, TX 77843, USA\\
and\\
Universit\'e Paris VI\\
Equipe d'Analyse, Case 186, 75252\\
Paris Cedex 05, France,\\ {pisier@math.tamu.edu}}

\date{}
\maketitle

\abstract{We prove that the 
operator Hilbert  space $OH$ does not embed
completely
 isomorphically  into 
the predual   of a semi-finite von~Neumann
algebra. This complements Junge's recent result
that it admits such an embedding in the non
semi-finite case. }

 \centerline{--------------------}
\bigskip

In remarkable recent work \cite{J}, Marius Junge
proves that the operator Hilbert  space $OH$
(from \cite{P1}, see also \cite{P2}) embeds completely
 isomorphically  (c.i.\ in short) into 
the predual $M_*$ of a von~Neumann algebra $M$ 
which is of type $III$; thus this algebra $M$ is
not semi-finite. In this note, we show that no
such embedding  can exist when $M$ is semi-finite.

The results we just stated all belong to the
currently very active field of ``operator spaces"
for which we refer the reader to the monographs
\cite{ER2,P3}. We merely recall a few basic facts
relevant for the present note. An operator space
is a Banach space given together with an isometric
embedding $E\subset B(H)$ into the algebra $B(H)$
of all bounded operators on
a Hilbert space $H$. Using this embedding,
we  equip the space
$M_n(E)$ (consisting of the
  $n\times n$ matrices with entries in $E$)
with the norm induced by
the space $M_n(B(H))$, naturally identified
isometrically with $B(H\oplus\dots \oplus H)$.
Let $F\subset B(K)$ be another operator space. In
operator space theory, the morphisms
are the completely bounded (c.b. in short) linear
maps: A linear map
$u\colon\ E\to F$ is called c.b. if the mappings
$u_n\colon\ M_n(E)\to M_n(F)$ defined
by $[x_{ij}] \to [u(x_{ij})]$ are uniformly
bounded when $n$ ranges over all integers $\ge
1$, and the cb-norm is defined
as
$\|u\|_{cb} =\sup_n \|u_n\|.$
The resulting normed space
of all c.b. maps $u\colon\ E\to F$
equipped with the cb-norm is denoted
by $CB(E,F)$.  If $u$ is
invertible with c.b. inverse, then $u$ is called
a complete isomorphism. For any operator space
$E\subset B(H)$, the Banach dual $E^*$ can be
equipped with a specific operator space
structure, say
$E^*\subset B(\H)$, for which the natural
identification $M_n(E^*)\simeq CB(E,M_n)$ is
isometric. On the other hand, the complex
conjugate
$\overline E$ can obviously be viewed as an
operator space using the canonical embeddings
$\overline E\subset \overline {B(H)}\simeq
B(\overline{H})$. Let $I$ be any set. In \cite{P1}, we
exhibited an operator space $E$ which is isometric
to $\ell_2(I)$ as a Banach space and such that
the canonical isometry (associated to the inner
product) $E\simeq \overline {E^*}$ is a complete
isometry. The latter operator space (which is
uniquely characterized by the preceding self-dual
property) is denoted by $OH(I)$ 
(or simply by $OH$ when $I=\NN$). We call it  the
operator Hilbert space.

In Banach space theory (and in commutative
harmonic analysis) the existence of
an isomorphic (actually isometric) embedding
of  $\ell_2$ (or  $\ell_2(I)$) in an $L_1$-space
plays a very important role, in connection with
the Khintchine and Grothendieck inequalities.
However, the non-commutative version of the
Khintchine inequality due to Lust-Piquard and
the author
(see \cite[Th. 8.4.1]{P2}),   when properly
interpreted,     leads to the  embedding  of a
different (Hilbertian but not OH) operator space
into
$L_1$, namely the space called $R+C$ in \cite{P2}.
This essentially implies that $OH$ does not embed
(c.i.) into a commutative $L_1$ space (see \cite{JO}
for details and more general results). Thus the
question (raised in \cite{P1}) whether
$OH$ itself embeds (c.i.) in a {\it
non-commutative}
$L_1$ space remained open, and was only recently
solved affirmatively by Junge \cite{J}. By a
non-commutative 
$L_1$-space, we mean the predual $M_*$ of a von
Neumann algebra $M$, equipped with the operator
space structure induced by the dual $M^*$.
Recall  that $M$ is called semi-finite
if it admits a  normal and faithful
trace which is also semi-finite, that is, 
although it is not necessarily finite, it 
admits sufficiently many elements on which it
is finite (see \cite{T} for more).
 
Our main theorem below shows that if $OH$ embeds
in $M_*$, then $M$ cannot be semi-finite.
This proves the need for a type III 
 algebra in Junge's work, and perhaps explains
the delay in resolving this embedding problem. 

 We will use the following easy consequence of
\cite[Cor.~3.4]{P3}.

\begin{lem}\label{lem1}
 Let $I$ be any set, let $E = OH(I)$ and let $M$ be a 
von~Neumann algebra. Consider a linear map $v\colon \ E\to M_*$.
\begin{itemize}
\item[(i)] If $v$ is completely bounded (c.b.\ in short), there is a normal 
state 
$\varphi$ on $M$ such that
\end{itemize}
\begin{equation}\label{eq1}
\forall~x\in M\qquad \|v^*x\|^2 \le
K^2(\varphi(xx^*)\varphi(x^*x))^{1/2}
\end{equation}
\begin{itemize}
\item[]
where $K=2^{9/4}\|v\|_{cb}$.
\item[(ii)] Conversely if there are $f_1,f_2$ in $M^*_+$ such that
\[
\forall~x\in M\hspace{1in} \|v^*x\|^2 \le
(f_1(xx^*)f_2(x^*x))^{1/2},
\]
 then necessarily
\end{itemize}
\begin{equation}\label{eq2}
\|v\|_{cb} \le (f_1(1)f_2(1))^{1/4}.
\end{equation}
\end{lem}

\begin{proof}
By \cite[Cor.~3.4]{PS} there is a state $f$ on
$M$ such that
\begin{equation}\label{eq3}
\forall~x\in M\hspace{1in}\|v^*x\|^2 \le K^2(f(xx^*)f(x^*x))^{1/2}.
\end{equation}
 We will use an argument
that can be traced back to \cite[p. 352]{R} and, in
the  non-commutative case, to \cite[Prop.~2.3]{H}.
Let $f = f_n+f_s$ be the decomposition  of $f$
into its normal and singular part. We set
$\varphi = f_n$. As explained  in the proof of
Prop.~2.3  in \cite{H}, there is an increasing net
$(p_\alpha)$ of  projections in $M$ such that
$p_\alpha\to 1$, say in the strong operator 
topology (SOT in short), and $f_s(p_\alpha) = 0$
for all $\alpha$. Note that (by  the
SOT-continuity of the product of
$M$ on bounded sets) we have $p_\alpha  xp_\alpha
\to x$, $p_\alpha xp_\alpha x^*p_\alpha\to xx^*$
and $p_\alpha  x^*p_\alpha xp_\alpha \to x^*x$
for the SOT.
For any $\xi$ in the unit ball of $E^*$, by assumption $x\to \langle\xi, 
v^*x\rangle$ is in $M_*$ (i.e.\ is ``normal''), hence
\begin{align*}
\langle \xi, v^*x\rangle &= \lim \langle \xi,
v^*(p_\alpha xp_\alpha)\rangle \\
\intertext{hence by \eqref{eq3}}
&\le K^2\lim(f(p_\alpha xp_\alpha x^*p_\alpha) f(p_\alpha x^*p_\alpha 
xp_\alpha))^{1/2}
\end{align*}
but since $f_s(p_\alpha) = 0$, a fortiori 
$f_s(p_\alpha xp_\alpha x^*p_\alpha) = 0=
f_s(p_\alpha x^* p_\alpha xp_\alpha)$, hence we
obtain
\[
|\langle\xi, v^*x\rangle| \le K^2\lim (\varphi(p_\alpha xp_\alpha x^*p_\alpha) 
\varphi(p_\alpha x^*p_\alpha xp_\alpha))^{1/2}
\]
hence we conclude
\[
|\langle\xi, v^*x\rangle| \le K^2(\varphi(xx^*)
\varphi(x^*x))^{1/2}
\]
 which immediately yields
(i). The proof of (ii) is identical to that of
the last  assertion in Corollary~3.4 in \cite{PS}.
\end{proof}

Our main result is the following.

\begin{thm}
Let $E = OH(I)$ with $I$
 an infinite set. Let $F\subset M_*$ 
be a subspace of  the predual of 
a semi-finite von~Neumann algebra $M$. Then for 
any c.b.\ maps
\[
u\colon \ E\to F \quad \text{and}\quad w\colon \ F\to E
\]
the composition $wv$ is compact.
\end{thm}

\begin{proof}
 Clearly this reduces to $I = \NN$ if we wish.
Let $\tau$ be
 a normal faithful semi-finite trace on $M$.

We will argue by contradiction. Assume that $wv$ is not compact. Then by the 
homogeneity of $OH$ (see \cite[p.~18]{P1}) we may assume that $wv$ is diagonal 
or 
even furthermore that $wv$ is the identity on $E$. In other words, we may as 
well assume $u$ invertible and $w=u^{-1}$.

Let $v\colon \ E\to M_*$ be the same map as 
$u$ but viewed as acting into $M_*$. 
By Lemma~\ref{lem1}, there is a normal state $\varphi$
 such that \eqref{eq1} holds. Let $e$ be the support
projection of $\varphi$ (i.e. we have 
$\varphi(1-e)=0$ and $\varphi(q)>0$ for any non
zero projection $q$ in $M$ with $q\le e$).
Then \eqref{eq1} implies that
 for any x in $M$, we have
$v^*(x(1-e))=0=v^*((1-e)x)$ hence
\[
v^*(x)=v^*(exe).
\]
 Thus if we replace $M$
by $eMe$ and  $\varphi$ by  $e\varphi $, we may
assume in addition that  $\varphi$ is faithful.

 Since 
$M_*\simeq L_1(\tau)$,
 we may assume $\varphi = \psi\cdot \tau$. Fix 
$0<\delta<1$.
 Let $p$ be the spectral projection associated to
$\psi$  with  respect to the set
$[\delta,\delta^{-1}]$ (for more details,
see e.g.\ \cite[p.~338]{SZ} or also \cite{N}).
 Note that  $\delta\tau(p) \le 
\tau(\psi) = \varphi(1) \le 1$, so 
that in particular $\tau(p)<\infty$
and moreover
\begin{equation}\label{eq4}
p\cdot \psi = 
\psi\cdot p = p\cdot \psi \cdot p 
\le \delta^{-1}
p.
\end{equation}
  On the 
other hand, let $\vp(\delta) =
\tau((1-p)\psi)=\varphi(1-p)$. Clearly (since
$\varphi$ is  faithful) we have
$\vp(\delta)\to 0$ when $\delta\to 0$. Thus if we
set, for all 
$y$ in $E$,
\[
v_\delta(y) = p\cdot v(y)\cdot p \text{ and } T_1(y) = v(y) (1-p), T_2(y)  =
(1-p)v(y)p,
\]
we have $v=v_\delta +
T_1+T_2$. We will show that $T_1+T_2$ is  small
when
$\delta\to 0$ so that $v_\delta$ can be viewed as
a perturbation of 
$v$. Indeed, for $x$ in $M$, we have
 $T^*_1(x) = v^*((1-p)x)$ and $x^*(1-p)x\le 
x^*x$ hence by \eqref{eq1} we have
\[
\|T^*_1(x)\|^2 \le K^2(\varphi((1-p)xx^*(1-p))
\varphi(x^*x)))^{1/2}
\]
 hence by \eqref{eq2}
\begin{align*}
\|T_1\|_{cb} &\le K(\varphi(1-p))^{1/4}\\
&\le K(\vp(\delta))^{1/4}.
\end{align*}
Similarly $\|T_2\|_{cb} \le
K(\vp(\delta))^{1/4}$, hence 
$\|v-v_\delta\|_{cb}\le \|T_1\|_{cb} + 
\|T_2\|_{cb} \le 2K(\vp(\delta))^{1/4}$. 
Let $\vp_0 = \|u^{-1}\|^{-1}_{cb}$. Clearly, if we choose $\delta$ small enough 
so that $2K(\vp(\delta))^{1/4}<\vp_0$ we have
$\|v-v_\delta\|_{cb} < \vp_0$,  hence, by an
elementary reasoning (based solely on the triangle
inequality in $M_n(M_*)$), the map
$v_\delta\colon
\ E\to pM_*p
\subset M_*$ is a completely isomorphic 
embedding. But now by \eqref{eq4} we have
\[
\|v^*_\delta (x)\|^2 =\|v^* (pxp)\|^2\le
K^2\delta^{-1}\|pxp\|_{L^2(\tau)} 
\|px^*p\|_{L^2(\tau)}
\]
hence since $\tau$ is {\it tracial} (this is 
where we make  crucial use of the semi-finiteness
assumption)
\[
\|v^*_\delta (x)\|^2 \le
K^2\delta^{-1}\|pxp\|^2_{L^2(\tau)}\le
K^2\delta^{-1}\min\{\|px\|^2_{L^2(\tau)},\|xp\|^2_{L^2(\tau)}  
 \}\le
K^2\delta^{-1} \|px\|^2_{L^2(\tau)}.
 \]
  By Lemma 2 below, since $\tau(p)<\infty$,
this is impossible.
\end{proof}

\begin{lem}\label{lem2}
 With the above notation.   Let $V\colon \ OH\to 
L_1(\tau)$ be a linear map for which there is $a$ in the unit ball of 
$L_2(\tau)$ and a  constant $C$ such that, for any $x$ in $M$, we have
\begin{equation}\label{eq5}
\|V^*(x)\|\le C\|ax\|_{L^2(\tau)}.
\end{equation}
Then,  for any isometry $J\colon \ C\to OH$,  $VJ$ is completely bounded 
from 
$C$ to $L_1(\tau)$.
In particular $V$ cannot be a completely isomorphic embedding.
\end{lem}

\begin{proof}
 By \eqref{eq5},  for any finite sequence 
$(x_i)$ 
in $M$
\begin{equation}\label{eq6}
\sum\|(VJ)^*(x_i)\|^2=\sum\|V^*(x_i)\|^2 \le C^2 \left\|\sum
x_ix^*_i\right\|.
\end{equation}
By a well known argument from \cite{ER1}, it follows
that
$
\|(VJ)^*\colon \ M\to C^*\|_{cb} \le C$, and hence
$VJ$ is c.b. .
Finally, if $V$ were a completely isomorphic
embedding, then
$VJ$ (when viewed as acting into the range of $V$)
would be   a
completely  bounded map from $C$ to $OH$, hence would
be (by Remark 2.11 in \cite{P1}) in the Schatten class $S_4$, and a
fortiori would be compact.
 But then $J$ itself would have to be compact which is
absurd.
\end{proof}

\n {\it Remarks.} 
\begin{itemize}
\item[(i)] Junge \cite{J} proves
that $OH_n$ embeds c.i. (with uniform constants)
into the predual of a finite dimensional
(hence semi-finite!) von Neumann algebra. More
precisely, he proves that there is $C>0$ so that
for any
$n$,  there is an integer $N$, a subspace
$F_n\subset M_N^*$  and a (complete)
isomorphism  $u_n\colon\ OH_n\to F_n$
 such that $\sup_n \|u_n\|_{cb}
\|u_n^{-1}\|_{cb}\le C$. It would be interesting
to   estimate $N$ as a function of $n$. 
\item[(ii)] The non-existence of
embeddings of $OH$ into $M_*$ when $M$ is
commutative is rather easy. In that case, even
the finite dimensional case (as in the preceding
point) is ruled out (see \cite{P4} for related
facts). The paper \cite{JO} contains stronger results
in the same direction.
\item[(iii)] The above
theorem remains valid with essentially the same proof for
$E=(R,C)_{\theta}$ (with $0< 
\theta<1$) in the sense of \cite{P2}, but this 
requires
  the   generalized version of Lemma
\ref{lem1} proved in \cite{P5}. This implies that, 
for any $1<p<2$, the
Schatten classes $S_p$ (and hence most
non-commutative $L_p$-spaces) do not embed
(c.i.)   into 
the predual   of any semi-finite von~Neumann
algebra.
\item[(iv)] Let $N_*$ be the predual of an injective factor
of type III and let $M$ be a (semi-finite)
von Neumann algebra of type ${\rm II}_\infty$.
 Junge proved that $OH$ embeds
completely isomorphically into $N_*$. 
Hence Theorem 1 implies that $N_*$
does not embed c.i. into $M_*$.
 This gives a somewhat partial answer
to the (still open) question raised in \cite{HRS}
of the existence of a (Banach space sense) isomorphic 
embedding of 
$N_*$ into $M_*$.

\end{itemize}

\end{document}